\documentstyle[12pt]{article}

\begin{document}
\newcommand{\ol }{\overline}
\newcommand{\ul }{\underline }
\newcommand{\ra }{\rightarrow }
\newcommand{\lra }{\longrightarrow }
\newcommand{\ga }{\gamma }
\newcommand{\st }{\stackrel }
\newcommand{\scr }{\scriptsize }
\title{{\Large\bf On the Order of Nilpotent Multipliers of Finite
$p$-Groups}\footnote{This research was in part supported by a
grant from IPM.}}
\author{ Behrooz Mashayekhy and Mahboobeh Alizadeh Sanati \\
Department of Mathematics, Ferdowsi University of Mashhad,\\
P.O.Box 1159-91775, Mashhad, Iran \\
Institute for Studies in Theoretical Physics and Mathematics, \\
Tehran, Iran}
\date{ }
\maketitle
\begin{abstract}
Let $G$ be a finite $p$-group of order $p^n$. YA. G. Berkovich
(Journal of Algebra {\bf 144}, 269-272 (1991)) proved that $G$ is
elementary abelian $p$-group if and only if the order of its
Schur multiplier, $M(G)$, is at the maximum case. In this paper,
first we find the upper bound $p^{\chi_{c+1}{(n)}}$ for the order
the $c$-nilpotent multiplier of $G$, $M^{(c)}(G)$, where
$\chi_{c+1}{(i)}$ is the number of basic commutators of weight
$c+1$ on $i$ letters. Second, we obtain the structure of $G$, in
abelian case, where $|M^{(c)}(G)|=p^{\chi_{c+1}{(n-t)}}$, for all
$0\leq t\leq n-1$. Finally, by putting a condition on the kernel
of the left natural map of the generalized Stallings-Stammbach
five term exact sequence, we show that an arbitrary finite
$p$-group with the $c$-nilpotent multiplier of maximum order is
an elementary abelian $p$-group.\\
{\it Key words and phrases}: Nilpotent multipliers, finite
$p$-group, elementary abelian $p$-group.\\ {\it
A.M.S.Classification 2000} : 20D15, 20E10, 20K01.
\end{abstract}
{\bf 1.Introduction and Motivation}

 Let $G$ be a group with a free presentation
$$ 1\lra R\lra F\lra G\lra 1\ ,$$
where $F$ is a free group and $R$ a normal subgroup of $F$ such that the above
sequence is exact. This paper focuses on the abelian group
$$ M^{(c)}(G)=\frac {R\cap \ga_{c+1}(F)}{[R, _cF]}\ \ ,$$
where $ \ga_{c+1}(F)$ is the $(c+1)$st term of the lower central series of $F$ and
$$ [R, _1F]=[R,F]\ ,$$
$$ [R, _cF]=[[R, _{c-1}F],F]\ .$$

 R. Baer [1] showed that $M^{(c)}(G)$ is independent of the choice of the free
presentation of $G$. The group $M^{(1)}(G)=M(G)$ is the much studied {\it Schur
multiplier} of $G$. We shall call $M^{(c)}(G)$ the $c$-nilpotent multiplier of
$G$. In fact, $M^{(c)}(G)$ is the well-known notion the {\it Baer invariant} of
$G$ with respect to the variety of nilpotent groups of class at most $c$, ${\cal
N}_c$. We also denote $M^{(c)}(G)$ by ${\cal N}_cM(G)$ (See [11]).

 In 1907, I. Schur [16], using representation method, showed that if $G$ is the
direct product of $A$ and $B$, then the Schur multiplier of
$G$ has the following form:
   $$ M(G)=M(A\times B)\cong M(A)\oplus M(B)\oplus (A_{ab}\otimes B_{ab}).$$
 A useful consequence of the above fact is to obtain an explicit formula for the
structure of the Schur multiplier of a finite abelian group, in particular for
an elementary abelian $p$-group.

 In 1956, J. A. Green [6] proved that if $G$ is a finite $p$-group of order
$p^n$, then $M(G)$, the Schur multiplier of $G$, is of order $p^{m(G)}$, where
$m(G)\leq\frac{1}{2}n(n-1)$.

 In 1991, YA. G. Berkovich [3], using the above facts and some properties of
covering groups proved that a finite $p$-group $G$ of order $p^n$
is elementary abelian $p$-group if and only if the order of
$M(G)$ is $p^{m(G)}$, where $m(G)=\frac{1}{2}n(n-1)$.

 Now, in this paper, we are going to concentrate on the Berkovich's
result for the $c$-nilpotent multiplier of a finite $p$-group.\\
 {\bf 2.Notation and Prelimineries}

 In this section, we are going to recall some notions which are used in
future.\\
{\bf Definition 2.1} [7]

 The notion of {\it basic commutators} on letters $x_1,x_2,\ldots , x_n,\ldots$
, are defined as follows:\\
$(i)$ The letters $x_1,x_2,\ldots ,x_n,\ldots $ are basic commutators of weight
one, ordered by setting $x_i<x_j$ if $i<j$.\\
$(ii)$ If basic commutators $c_i$ of weight $wt(c_i)<k$ are defined and
ordered, then define basic commutators of weight $k$ by the following rules:\\
$[c_i,c_j]$ is a basic commutator of weight $k$ if

 1. $wt(c_i)+wt(c_j)=k\ ,$

 2. $c_i>c_j\ ,$

 3. if $c_i=[c_s,c_t]$, then $c_j\geq c_t$ .\\
Then continue the order by setting $c>c_i$ whenever $wt(c)>wt(c_i)$ and
fixing any order amongst those of weight $k$ and finally numbering them in
order.\\
{\bf Theorem 2.2} (P.Hall [7,8])

 Let $F=<x_1,x_2,\ldots ,x_d>$ be a free group, then
$$ \frac {\ga_n(F)}{\ga_{n+i}(F)} \ \ , \ \ \ \  1\leq i\leq n$$
is the free abelian group freely generated by the basic commutators of weights
$n,n+1,\ldots ,n+i-1$ on the letters $\{x_1,\ldots ,x_d\}.$\\
{\bf Theorem 2.3} (Witt Formula [7])

 The number of basic commutators of weight $n$ on $d$ generators is  given
by the following formula:
$$ \chi_n(d)=\frac {1}{n} \sum_{m|n}^{} \mu (m)d^{n/m}$$
where $\mu (m)$ is the {\it Mobious function}, and defined to be
   \[ \mu (m)=\left \{ \begin{array}{ll}
      1 & {\rm if}\ m=1, \\ 0 & {\rm if}\ m=p_1^{\alpha_1}\ldots
p_k^{\alpha_k}\ \ \exists \alpha_i>1, \\ (-1)^s & {\rm if}\ m=p_1\ldots p_s,
\end{array} \right.  \]
{\bf Definition 2.4}

 Let $G$ be a group with a free presentation
             $$ 1\lra R\lra F\lra G\lra 1\ , $$
where $F$ is a free group. Then the {\it Baer invariant} of $G$ with respect
to the variety $\cal V$, denoted by ${\cal V}M(G)$ is defined to be
$$ {\cal V}M(G)=\frac {R\cap V(F)}{[RV^*F]}\ . $$

 One can see that the Baer invariant of the group $G$ is always abelian and
independent of the choice of the presentation of $G$.
( See C. R. Leedham-Green and S. McKay [11], from where we have taken our
notation and H. Neumann [15] for the notion of varieties of groups.)

 In particular, if $\cal V$ is the variety of abelian groups, $\cal A$, then the
Baer invariant of the group $G$ will be
  $$ \frac {R\cap F'}{[R,F]}\ , $$
which is isomorphic to {\it the Schur multiplier} of $G$,
denoted by $M(G)$, where, in finite case, is the {\it second cohomology group}
of $G$, $H^2(G,{\bf C}^*)$. Also if $\cal V$ is the
variety of nilpotent groups of class at most $c\geq 1$, ${\cal N}_c$, then the
Baer invariant of the group $G$ will be
  $$ {\cal N}_cM(G)=\frac {R\cap \ga_{c+1}(F)}{[R,\ _cF]}\ .$$
We also call it the $c$-nilpotent multiplier of $G$.\\
{\bf Definition 2.5}

 A variety $\cal V$ is said to be a {\it Schur-Baer variety} if for any group
$G$ for which the marginal factor group, $G/V^*(G)$, is finite, then the verbal
subgroup,$V(G)$, is also finite and $|V(G)|$ divides a power of $|G/V^*(G)|$.
I.\ Schur in [16] proved that the variety of abelian groups, $\cal A$, is a
Schur-Baer variety. Also R. Baer in [2] proved that the variety defined by some
outer commutator words has the above property.
 The following theorem tells us a very important property of Schur-Baer
varieties.\\
{\bf Theorem 2.6}(C. R. Leedham-Green and S. McKay [11])

 The following conditions on the variety $\cal V$ are equivalent:\\
$(i)$ $\cal V$ is a Schur-Baer variety.\\
$(ii)$ For every finite group $G$, its Baer invariant, ${\cal V}M(G)$, is of
order dividing a power of $|G|$\ .\\
{\bf Theorem 2.7} (I. Schur [10, 16])

 Let
$$ G\cong {\bf Z}_{n_1}\oplus {\bf Z}_{n_2}\oplus \ldots \oplus {\bf Z}_{n_k}$$
be a finite abelian group,
where $n_{i+1}|n_i$ for all $1\leq i\leq k-1$ and $k\geq 2$. Then
$$ M(G)\cong {\bf Z}_{n_2}\oplus {\bf
Z}_{n_3}^{(2)}\oplus \ldots \oplus {\bf
Z}_{n_k}^{(k-1)} \ ,$$
where ${\bf Z}_n^{(m)}$ denote the direct sum of $m$ copies of ${\bf Z}_n$ .\\
{\bf Theorem 2.8} (M. R. Jones [9, 10])

 Let $G$ be a finite $p$-group of order $p^n$ and let $G$ admit a presentation
with $d$ generators. Then
$$ p^{1/2d(d-1)}\leq |M(G)||G'|\leq p^{1/2n(n-1)}\ .$$
{\bf Theorem 2.9} (J. A. Green [6])

 Let $G$ be a finite $p$-group of order $p^n$. Then the order of $M(G)$ is at
most $p^{1/2n(n-1)}$ .\\
{\bf Theorem 2.10} (YA. G. Berkovich [3])

 Let $G$ be a finite $p$-group of order $p^n$. Then $|M(G)|=p^{1/2n(n-1)}$ if and
only if $G$ is an elementary abelian $p$-group.\\
{\bf Definition 2.11} [12]

 By a {\it pre-crossed} module we mean a pair of groups $(A,G)$ and
homomorphism
$\alpha :A\lra G$, together with an action of $G$ on $A$ which satisfies the
condition
$$ \alpha(a^g)=g^{-1}\alpha(a)g\ \ \ for\ all\ g\in G\ and\ a\in A\ .$$
A pre-crossed module for which, in addition, we have
$$ a^{\alpha(b)}=b^{-1}ab\ \ \ for\ all\ a,b\in A $$
is called a {\it crossed module}.

 Let $A$ be a $G$-precrossed module, then the non-abelian tensor product
$A\otimes G$ is defined to be the group with generators
$$ \{ <a,x>|a\in A,\ x\in G\}$$
subject to the relations
$$ <a_1a_2, x>=<a_2^{\alpha(a_1)},x^{\alpha(a_1)}><a_1,x> $$
$$<a,x_1x_2>=<a,x_1><a^{x_1},x_2^{x_1}>\ \ .$$
\hspace{-.2in}{\bf Theorem 2.12} (J. Stallings and U. Stammbach [4])

 Every extension $1\ra N\ra G\ra Q\ra 1$ determines an exact sequence
$$M(G)\lra M(Q)\lra N/[N,G]\lra G/G'\lra Q/Q'\lra 1\ \ .$$
{\bf 3. The Main Results}

 In order to present some main results we need the following theorems.\\
{\bf Theorem 3.1} (B. Mashayekhy and M. R. R. Moghaddam [13])

 Let $ G\cong {\bf Z}_{n_1}\oplus {\bf Z}_{n_2}\oplus \ldots \oplus {\bf
Z}_{n_k}$ be a finite abelian group,
where $n_{i+1}|n_i$ for all $1\leq i\leq k-1$ and $k\geq 2$.
Then, for all $c\geq 1$, the $c$-nilpotent multiplier of $G$ is
$$ {\cal N}_cM(G)=M^{(c)}(G)\cong {\bf Z}_{n_2}^{(b_2)}\oplus {\bf
Z}_{n_3}^{(b_3-b_2)}\oplus \ldots \oplus {\bf
Z}_{n_k}^{(b_k-b_{k-1})} \
\ ,$$ where $b_i$ is the number of basic commutators of weight $c+1$ on
$i$ letters. Note that ${\bf Z}_n^{(m)}$ denotes the direct sum of $m$ copies of
${\bf Z}_n$.\\
{\bf Theorem 3.2} (M. R. R. Moghaddam [14])

 Let $\cal V$ be a variety of polynilpotent groups and let $G$ be a $p$-group of
order $p^n$, where $p$ is any prime number. Then
$$ |{\cal V}M(G)||V(G)|\leq |{\cal V}M({\bf Z}_p^{(n)})|\ \ .$$

 Now, we can get an upper bound for the order of a $c$-nilpotent multiplier of a
finite $p$-group as follows.\\
{\bf Corollary 3.3}

 Let $G$  be a $p$-group of order $p^n$. Then the order of $M^{(c)}(G)$ is at
most $p^{\chi_{c+1}(n)}$, where $\chi_{c+1}(n)$ is the number of basic
commutators of weight $c+1$ on $n$ generators.\\
\hspace{-.2in}{\bf Proof.}

 Since
$ {\bf Z}_p^{(n)}\cong {\bf Z}_p\oplus {\bf Z}_p\oplus \ldots \oplus {\bf Z}_p$
 (n-copies), by Theorem 3.1 we have
$$ M^{(c)}({\bf Z}_p^{(n)})\cong {\bf Z}_p^{(b_2)}\oplus {\bf
Z}_p^{(b_3-b_2)}\oplus \ldots \oplus {\bf Z}_p^{(b_k-b_{k-1})}\ \ . $$
Therefore $|M^{(c)}({\bf Z}_p^{(n)})|=p^{b_2+(b_3-b_2)+\ldots
+(b_n-b_{n-1})}=p^{b_n}=p^{\chi_{c+1}(n)}$.

 Now, by Theorem 3.2 we have
$$ |M^{(c)}(G)|\leq |M^{(c)}(G)||\ga_{c+1}(G)|\leq |M^{(c)}({\bf
Z}_p^{(n)})|=p^{\chi_{c+1}(n)}\ .\Box$$

 Note that Corollary 3.3 is a vast generalization of Green's result (Corollary
2.9), since
$\chi_2(n)=1/2\sum_{m|2}^{}\mu(m)n^{2/m}=1/2(\mu(1)n^2+\mu(2)n^1)=1/2(n^2-n)=1/2
n(n-1)$.

 Now in the following we are going to present a vast generalization of the Berkovich's result
(Theorem 2.10) in abelian case.\\
{\bf Theorem 3.4}

Let $G$ be an abelian $p$-group of order $p^n$. Then for all $c\geq 1$ and all $0\leq t\leq n-1$
\begin {center}
$|M^{(c)}(G)|=p^{\chi_{c+1}(n-t)}$ if and only if $G\cong {\bf Z}_{p^{t+1}}\oplus\underbrace{ {\bf Z}_p\oplus \ldots\oplus {\bf Z}_p}_{n-t-1-copies}\ \ .$
\end{center}
{\bf Proof.}

Let $G\cong{\bf Z}_{p^{t+1}}\oplus\underbrace{{\bf
Z}_p\oplus\ldots\oplus{\bf Z}_p}_{n-t-1-copies}$\ . Then, by
Theorem 3.1, we have $M^{(c)}(G)\cong{\bf Z}_p^{(b_2)}\oplus{\bf
Z}_p^{(b_3-b_2)}\oplus\ldots\oplus{\bf
Z}_p^{(b_{n-t}-b_{n-t-1})}\ .$ Hence
$$|M^{(c)}(G)|=p^{b_2+(b_3-b_2)+\ldots+(b_{n-t}-b_{n-t-1})}=p^{b_{n-t}}=p^{\chi_{c+1}(n-t)}\ .$$
 Conversely, suppose that $|M^{(c)}(G)|=p^{\chi_{c+1}(n-t)}\ .$ Since $G$ is an
abelian $p$-group of order $p^n$, so by the fundamental theorem of finitely
generated abelian groups we have
$$G\cong{\bf Z}_{p^{\alpha_1}}\oplus{\bf Z}_{p^{\alpha_2}}\oplus\ldots\oplus{\bf Z}_{p^{\alpha_k}}\ \ ,$$
where $\alpha_1\geq\alpha_2\geq\ldots\geq\alpha_k$ and $\alpha_1+\alpha_2+\ldots+\alpha_k=n$ . By Theorem 3.1 $M^{(c)}(G)\cong{\bf Z}_{p^{\alpha_2}}^{(b_2)}\oplus{\bf Z}_{p^{\alpha_3}}^{(b_3-b_2)}\oplus\ldots\oplus{\bf Z}_{p^{\alpha_k}}^{(b_k-b_{k-1})}\ .$ Thus
$$|M^{(c)}(G)|=p^{\alpha_2b_2+\alpha_3(b_3-b_2)+\ldots+\alpha_k(b_k-b_{k-1})}\ \ .$$
 On the other hand, by hypothesis $|M^{(c)}(G)|=p^{\chi_{c+1}(n-t)}=p^{b_{n-t}}$. Therefore we have $ b_{n-t}=\alpha_2b_2+\alpha_3(b_3-b_2)+\ldots +\alpha_k(b_k-b_{k-1})\ $. Since $\alpha_2\geq\alpha_i$ for all $2\leq i\leq k$, we have
$$b_{n-t}\leq\alpha_2b_2+\alpha_2(b_3-b_2)+\ldots+\alpha_2(b_k-b_{k-1})=\alpha_2b_k\hspace{1cm}(I)\ .$$
Also $\alpha_k\leq\alpha_i$ for all $2\leq i\leq k$, so we have
$$b_{n-t}\geq\alpha_kb_2+\alpha_k(b_3-b_2)+\ldots+\alpha_k(b_k-b_{k-1})=\alpha_kb_k\hspace{1cm}(II)\ .$$
If $n=1$ then $t=0$. So we have $G\cong{\bf Z}_p$, $|M^{(c)}(G)|=1=p^{b_{1}}$. Hence the result holds\ .\\
If $n=2$ then there are two cases $t=0$ or $t=1$. If $t=0$, then $|M^{(c)}(G)|=p^{b_2}$ and so $G\cong{\bf Z}_p\oplus{\bf Z}_p$. If $t=1$ then $|M^{(c)}(G)|=p^{b_1}=1$ and $G\cong{\bf Z}_{p^2}$. Hence the result holds. Therefore we can assume that $n\geq 3$.

Now, we claim that $k=n-t$. In order to prove the claim, consider the following two cases.

Case one: If $k>n-t$, then clearly we have the inequalities $\alpha_kb_k\geq b_k>b_{n-t}$ which contradicts to inequality $(II)$.

Case two: If $k<n-t$. Put $k=n-t-j$ for some $1\leq j<n-t-1$. Clearly $\alpha_2$ is maximum when $\alpha_1=\alpha_2$ and $\alpha_3=\ldots=\alpha_k=1\ .$ So in this case $2\alpha_2+k-2=n$. Hence we have $\alpha_2\leq(t+j+2)/2\ .$ So by $(I)$ we have the following inequality:
$$b_{n-t}\leq\frac{t+j+2}{2}b_{n-t-j}\ .$$
Now, we need the following inequality which is proved by induction on $n\geq 3$.
$$t+j+2\leq 2(n-(t+1))\ldots(n-(t+j))\hspace{2cm} (III)\ ,$$
where $0\leq t\leq n-1 and 1\leq j<n-t-1 .$ If $n=3$, then $t=0$ and $j=1$. So $0+1+2\leq 2(3-1)=4$ .

Now, let $(III)$ holds for $n$ with the above condition and we investigate the following inequality:
$$t+j+2\leq 2(n-t)\ldots(n+1-(t+j)),\ \ 0\leq t<n,\ \ 1\leq j<n-t\ .$$
 If $0\leq t<n-1$ and $1\leq j<n-t-1$, then by induction hypothesis we have
$$t+j+2\leq 2(n-(t+1))\ldots(n-(t+j))\leq 2(n-t)\ldots(n+1-(t+j))\ .$$
 If $0\leq t<n-1$ and $j=n-t-1$, then $k=1$ and so $G\cong{\bf Z}_{p^n}$. Hence $|M^{(c)}(G)|=1=p^{b_{n-t}}$. Thus $t=n-1$ and $j=0$ which is a contradiction to $j\geq 1$.

Also, if $t=n-1$ and $1\leq j<n-t$, then $1\leq j<n-(n-1)=1$, which is a contradiction. Hence the last two cases do not happen. Thus, the following inequalities hold by $(III)$:
$$b_{n-t}\leq\frac{t+j+2}{2}b_{n-t-j}\leq (n-(t+1))\ldots(n-(t+j))b_{n-t-j} \ .$$
By a routine combinatorial disscusion and considering the definition of basic
commutators we can see that $ib_i<b_{i+1}\ .$ (Note that we can choose $i$
letters from $i+1$ letters in $i$ cases.) Therefore we have
$$b_{n-t}\leq(n-(t+1))\ldots(n-(t+j))b_{n-t-j}<(n-(t+1))\ldots(n-(t+j-1))b_{n-(t+j-1)}$$
$$<\ldots <(n-(t+1))b_{n-(t+1)}< b_{n-t} \ .$$
which is a contradiction. Hence $k=n-t$ and our claim is true. So we can consider
$$G\cong{\bf Z}_{p^{\alpha_1}}\oplus{\bf Z}_{p^{\alpha_2}}\oplus\ldots\oplus{\bf Z}_{p^{\alpha_{n-t}}}\ .$$
 Now we show that $\alpha_1=t+1$ and $\alpha_2=\ldots=\alpha_{n-t}=1\ .$ We have $|M^{(c)}(G)|=p^{\alpha_2b_2+\alpha_3(b_3-b_2)+\ldots+\alpha_{n-t}(b_{n-t}-b_{n-t-1})}\ .$ By $(II)$ we have $\alpha_{n-t}b_{n-t}\leq b_{n-t}$ and so $\alpha_{n-t}=1$. Also, $\alpha_{n-t-1}\leq \alpha_i$ for $2\leq i\leq n-t-1$, thus
$$\alpha_{n-t-1}b_2+\alpha_{n-t-1}(b_3-b_2)+\ldots+\alpha_{n-t-1}(b_{n-t-1}-b_{n-t-2})+(b_{n-t}-b_{n-t-1})\leq b_{n-t}\ .$$
Hence we can conclude that $\alpha_{n-t-1}=1$. Similarly, one can see that  $\alpha_2=\ldots=\alpha_{n-t-1}=1.$ Since $\alpha_1+\ldots+\alpha_{n-t}=n,$ we have $\alpha_1=t+1$, as required. $\Box$\\

Now, the following corollary is a generalization of Berkovich's result (Theorem 2.10).\\
{\bf Corollary 3.5}

 Let $G$ be an abelian $p$-group of order $p^n$. Then
$|M^{(c)}(G)|=p^{\chi_{c+1}(n)}$ if and only if $G$ is the elementary abelian
$p$-group.\\

Now, in order to deal with non-abelian case, we need the following two
important results.\\
{\bf Theorem 3.6} (A. Fr\"{o}hlich [5])

 Let $1\ra N\ra G\ra Q\ra 1$ be a $\cal V$-central extension, where $\cal V$ is
any variety of groups (this means that the above sequence is exact and $N$ is
contained in the marginal subgroup of $G$, $V^*(G)$). Then the following
five-term exact sequence exists:
$$ {\cal V}M(G)\st {\theta}{\lra} {\cal V}M(Q)\lra N\lra G/V(G)\lra Q/V(Q)\lra
1\ \ .$$
{\bf Theorem 3.7} (A. S.-T. Lue [12])

 Let $1\ra N\ra G\ra Q\ra 1$ be an ${\cal N}_c$-central extension. Then the
following sequence is exact:
$$N\otimes \underbrace{G/\ga_{c+1}(G)\otimes \ldots \otimes
G/\ga_{c+1}(G)}_{c-copies}\lra {\cal
N}_cM(G)\st {\theta}{\lra }{\cal N}_cM(Q)$$ $$\lra N\lra G/\ga_{c+1}(G)\lra
Q/\ga_{c+1}(Q)\lra 1\ \ ,$$
where the above tensor product is the non-abelian tensor product. Moreover
 $$ ker\theta =\frac {[S, _cF]}{[R, _cF]}\ \ ,$$
where $G=F/R\ ,\ Q=F/S$ are free presentations for $G$ and $Q$, respectively.\\
{\bf Theorem 3.8}

 Let $G$ be a finite $p$-group of order $p^n$. If
$|M^{(c)}(G)|=p^{\chi_{c+1}(n)}$, then\\
$(i)$ There is an epimorphism $M^{(c)}(G)\st {\theta}{\lra }M^{(c)}(G/G')$ which
is obtain from the Fr\"{o}hlich's sequence.\\
$(ii)$ If $ker\theta =1$, then $G$ is an elementary abelian $p$-group.\\
{\bf Proof.}

 $(i)$ By Theorem 3.2 and Corollary 3.3
$$ |M^{(c)}(G)||\ga_{c+1}(G)|\leq |M^{(c)}({\bf Z}_p^{(n)})|=p^{\chi_{c+1}(n)}\
\ .$$
Since $|M^{(c)}(G)|=p^{\chi_{c+1}(n)}$, we have $|\ga_{c+1}(G)|=1$.\\
Now, set $N=G'$ and consider the exact sequence
$$ 1\lra G'\lra G\lra G/G'\lra 1 \ \ .$$
Since $\ga_{c+1}(G)=1$, the above exact sequence is an ${\cal N}_c$-central
extension. Therefore by Theorem 3.6 we have the following exact sequence:
$$M^{(c)}(G)\st {\theta}{\lra }M^{(c)}(G/G')\st {\beta}{\lra}
G'\st{\alpha}{\lra} G\lra G/G'\lra 1\ \ .$$
Clearly $\alpha$ is monomorphism and so $Im\beta =1$. This means that $\theta$
is epimorphism.

$(ii)$ Let $ker\theta =1$, then
$$|M^{(c)}(G/G')|=|M^{(c)}(G)|=p^{\chi_{c+1}(n)}\ .$$
Since $|G|=p^n$, so $|G/G'|\leq p^n$. Hence by Corollary 3.3 we have
$|G/G'|=p^n$. This means that $G'=1$ and so $G$ is abelian. Now by Corollary 3.5
$G$ is an elementary abelian $p$-group. $\Box$

E-mail: mashaf@math.um.ac.ir
\end{document}